\def\lastrevised{May 26, 2011.}

\documentclass[12pt]{article}
\usepackage{amssymb}

\def\header{\today}
\def\header{A.Miller\hfill Maximum Principle \hfill}
\def\al{\alpha}
\def\be{\beta}
\def\checkcheck#1{\check{\check{#1}}}
\def\dom(#1){{\rm dom}(#1)}
\def\ff{{\mathcal F}}
\def\forces{{\;\Vdash}}
\def\ga{\gamma}
\def\gg{{\mathcal G}}
\def\ka{{\kappa}}
\def\la{\langle}
\def\name#1{\stackrel{\circ}{#1}}
\def\namecheck#1{\name{\check{#1}}}

\def\nn{{\mathcal N}}
\def\om{\omega}
\def\one{{\bf 1}}
\def\pnames{(\poset{-{\rm names}})^M}
\def\poset{{\mathbb P}}
\def\power(#1){{\mathcal P}(#1)}
\def\proof{\par\noindent Proof\par\noindent}
\def\pr{\prime}
\def\qed{\par\noindent QED\par\bigskip}
\def\ra{\rangle}
\def\res{\upharpoonright}
\def\rmand{\mbox{ and }}
\def\rmiff{{\mbox{ iff }}}

\def\si{\sigma}

\def\sqleq{\sqsubseteq}
\def\st{\;:\;} 
\def\su{\subseteq}
\def\tlG{\tilde{G}}
\def\tleq{\unlhd}
\def\tlpi{\tilde{\pi}}
\def\zz{{\mathbb Z}}

\newtheorem{theorem}{Theorem}
\newtheorem{example}[theorem]{Example}

\newtheorem{lemma}[theorem]{Lemma}
\newtheorem{prop}[theorem]{Proposition}
\newtheorem{define}[theorem]{Definition}

\begin{document}

\begin{center}
{\large The maximum principle in forcing and the axiom of choice}
\end{center}

\begin{flushright}
Arnold W. Miller
\footnote{
\par Mathematics Subject Classification 2000: 03E25 03E40
\par Keywords: Forcing, Maximum Principle, Maximal Antichains,
Countable Axiom of Choice, Dedekind finite
\par Last revised \lastrevised
}
\end{flushright}

\def\address{\begin{flushleft}
Arnold W. Miller \\
miller@math.wisc.edu \\
http://www.math.wisc.edu/$\sim$miller\\
University of Wisconsin-Madison \\
Department of Mathematics, Van Vleck Hall \\
480 Lincoln Drive \\
Madison, Wisconsin 53706-1388 \\
\end{flushleft}}

\begin{center}  Abstract  \end{center}
\begin{quote}
In this paper we prove that the maximum principle in forcing
is equivalent to the axiom of choice.  We also look
at some specific partial orders in the basic Cohen model.
\end{quote}

Lately we have been thinking about forcing over models
of set theory which do not satisfy the axiom of choice
(see Miller \cite{longbor,ded}).
One of the first uses of the axiom of choice in forcing
is:

\begin{center}
Maximum Principle \\
$p\forces \;\exists x \;\theta(x)\;\;\;\;$  iff
$\;\;\;$ there exists a name $\tau$
$\;\;\;p\forces \;\theta(\tau)$.
\end{center}

Recall some definitions.
For a partial order $\poset=(\poset,\tleq)$ and $p,q\in\poset$
we say that $p$ and $q$ are compatible iff there exists an $r\in\poset$
with $r\tleq p$ and $r\tleq q$.  Otherwise $p$ and $q$ are incompatible.
A subset $A\su\poset$ is an antichain iff any two distinct elements
of $A$ are incompatible.   It is maximal iff every $p\in\poset$
is compatible with some $q\in A$.

The standard definition of $p\forces \;\exists x \;\theta(x)$ is given by:
$$p\forces \;\exists x \;\theta(x)\;\;\; \rmiff \;\;\;
\forall q\tleq p\;\exists r\tleq q\;\exists \tau \;\; r\forces \theta(\tau)$$
here $p,q,r$ range over $\poset$ and $\tau$ is a $\poset$-name.
The usual proof of the maximum principle is to choose a maximal antichain $A$
beneath $p$ of such $r$ and then choose names $(\tau_r\st r\in A)$
such that $r\forces\theta(\tau_r)$ for each $r\in A$.
Finally name $\tau$ is constructed from $(\tau_r\st r\in A)$
in an argument which does not use the axiom of choice.
For details the reader is referred to Kunen \cite{kunen} page 226,
who calls it the Maximal Principle.

Shelah \cite{shelah} and Bartosyznski-Judah \cite{barto} refer to
the maximum principle as the ``Existential Completeness Lemma''.
Takeuti-Zaring \cite{takeuti} use ``Maximum Principle''
to title their Chapter 16.

Jech \cite{jech} uses boolean valued models to do
forcing proofs.  He refers to the boolean algebra version of
the maximum principle as: ``$V^B$ is full'' , see Lemma 14.19
p.211.  He notes that this is the only place in his
chapter where the axiom of choice is used.  

We don't know if anyone has ever wondered if the axiom
of choice is necessary to prove the maximum principle.
First note that the axiom of choice is needed to give
the first step of the proof: Finding a maximal antichain.

\begin{theorem} \label{equiv}
The axiom of choice is equivalent to the
statement that every partial order contains a maximal
antichain.
\end{theorem}
\proof
Let $(X_i\st i\in I)$ be any family of nonempty pairwise
disjoint sets.
Let $$\poset=\bigcup_{i\in I}\;\;\om\times X_i$$ strictly ordered by:
$(n,x)\lhd (m,y)$ iff $n>m$ and $\exists i\in I \;\;x,y\in X_i$.

Note that any maximal antichain must
consist of picking exactly one element out of each $\om\times X_i$.
Hence we get a choice function.
\qed

The partial order used here is trivial in the forcing sense.
What happens if we only consider partial orders in which
every condition has at least two incompatible extensions?

In the literature on the axiom of choice
there is a property called the Antichain Property (A).
However, it is antichain in the sense of pairwise
incomparable not pairwise incompatible.  The property
(A) states that every partial order contains a maximal subset
$A$ of pairwise incomparable elements (i.e. for all $p,q\in A$
if $p\tleq q$, then $p=q$).

In ZF property (A) is equivalent to the axiom
of choice (but unlike Theorem \ref{equiv}) property (A) is
strictly weaker in set theory with atoms, i.e.,
it holds in some Fraenkel-Mostowski permutation model
in which the axiom of choice is false.
These two results are due to H.Rubin \cite{rubin} and Felgner-Jech
\cite{felgnerjech}.  See Chapter 9 of Jech \cite{jechac}.

\begin{theorem} \label{failure}
The axiom of choice is equivalent to the maximum principle.
\end{theorem}
\proof
Let $(X_i\st i\in I)$ be any family of nonempty pairwise
disjoint sets.
Let $\poset=I\cup\{\one\}$ strictly ordered by
$i\lhd \one$ for each $i\in I$ and the elements of $I$ pairwise
incomparable.
As usual the standard names for elements of the ground
model are defined by induction
$$\check{x}=\{(\one,\check{y})\st y\in x\}$$
and
$$\name{G}=\{(p,\check{p})\st p\in\poset\}$$ is a
name for the generic filter.

Then
$$\one\forces \exists x (\;\exists i\in \check{I}\cap \name{G}
\;\;x\in \check{X}_i)$$
which we may write as:
$$\one\forces \exists x \theta(x).$$
Applying the maximum principle, there exists
$\poset$-name $\tau$ such that
$$\one\forces \theta(\tau).$$
Then for each $i\in I$ we would have
to have a unique $x_i\in X_i$ such that
$$i\forces\tau=\check{x_i}.$$
This gives us a choice function.
\qed

This partial order is also trivial from the forcing point of view.
A nontrivial partial order which works is
$$\poset= (I\times 2^{<\om}) \cup \{\one\}$$
which is forcing equivalent to $2^{<\om}$.  In either
of these examples one can show (without using the axiom of choice)
that every
dense subset contains a maximal antichain.  Hence we can
think of them as showing that the second use of the axiom
of choice in the proof of the maximum principle, the choosing of names,
is also equivalent to the axiom of choice.

Note that the maximum principle holds for the suborder
$I\su \poset$.  So the maximum principle could fail for a partial
order but hold for a dense suborder.

What can be proved without the axiom of choice in the ground model?
For example,
if a partial order can be well-ordered in type $\ka$ and
choice holds for families of size $\ka$, then the usual
proof of the maximal principle goes thru.

We note a special case for which the maximum principle holds.

\begin{prop}\label{wellorder}(ZF)
Suppose $\ka$ is an ordinal and
$$p\forces \exists \al<\check{\ka} \;\;\; \theta(\al)$$
then there exists a name $\tau$ such that
$$p\forces \theta(\tau)$$
\end{prop}
\proof
Take $\tau$ to be a name for the least ordinal satisfying $\theta$:
$$\tau=\{(q,\check{\be})\st q\tleq p\rmand \forall \ga\leq\be
\;\; q\forces \;\neg\;\theta(\check{\ga})\}.$$
\qed

\section*{Basic Cohen model}

The Basic Cohen model $\nn$ for the negation of the axiom of choice
is described in Cohen \cite{cohen} and Jech \cite{jechac}.
It is the analogue of Fraenkel's 1922 permutation model.

\bigskip
\noindent 
One could\footnote{Since this model is the original and
simplest model in which the axiom of choice fails,
we think it is interesting to study its
properties just for its own sake.} ask:
In $\nn$ which partial orders have the maximum principle?

\bigskip

\def\inj(#1,#2){Inj(#1,#2)}

\begin{define}
Given infinite sets $I$ and $J$ let
$\inj(I,J)$ be the partial order of finite injective maps
from $I$ to $J$, i.e., $r\in\inj(I,J)$ iff
$r\su I\times J$ is finite and $u,v)$
It is ordered by reverse inclusion: $r_1\tleq r_2$ iff
$r_1\supseteq r_2$.
\end{define}

Recall that in $\nn$ the failure of the countable
axiom of choice is witnessed by
an infinite Dedekind finite $X\su \power(\om)$.
We consider the following three partial orders:
$\inj(\om,\om)$, $\inj(X,X)$, and $\inj(\om,X)$.

We show that the maximum principle holds for one of these
partial orders and fails for the other two.
The easiest case is $\inj(\om,\om)$.  The following lemma takes
care of it.

\begin{lemma} Suppose that the countable axiom of choice fails and
$\poset$ is a nontrivial partial order which can be well-ordered.
Then $\poset$ fails to satisfy the maximum principle.
\end{lemma}
\proof
By nontrivial we mean that every condition has at least
two incompatible extensions.
Hence we can find $\la p_n\in\poset\st n\in\om\ra$
such that $p_n$ and $p_m$ are incompatible whenever
$n\neq m$.
Suppose $\{X_n\st n\in \om\}$ is a family of nonempty sets
without a choice function.
Note that
$$\one\forces \exists x \;\forall n\in\check{\om}\;\;
 (\check{p}_n \in\name{G} \to x\in \check{X}_n).$$
We claim that this is a witness for the failure of
the maximum principle.  Suppose not and
let $\tau$ be $\poset$-name for which
$$\one\forces \forall n\in\check{\om}\;\;
 (\check{p}_n \in\name{G} \to \tau\in \check{X}_n).$$

Since $\poset$ can be well-ordered, we may choose for
each $n$ a $q_n\tleq p_n$ and $x_n\in X_n$ such
that
$$q_n\forces \tau=\check{x}_n.$$
But this would give a choice function for the family
$\{X_n\st n\in \om\}$.
\qed

\begin{theorem}
In $\nn$ the maximum principle fails for $\inj(\om,\om)$.
\end{theorem}
\proof
This follows from the Lemma, since $\inj(\om,\om)$ is well-orderable
and nontrivial,
and the countable axiom of choice fails in $\nn$.
\qed
Of course, there are many partial orders for which this
applies.  We choose to highlight $\inj(\om,\om)$ because
it is simple and superficially similar to 
the other two partial orders $\poset_0=\inj(X,X)$ and 
$\poset_1=\inj(\om,X)$.

\begin{theorem}
In $\nn$ the maximum principle fails for $\poset_0=\inj(X,X)$.
\end{theorem}
\proof

We start with a description of $\nn$. Fix $M$ a countable
standard transitive model of ZFC.

Working in $M$ let
$\poset=Fn(\om\times\om,2,\om)$ be the poset of
finite partial functions, i.e., $p\in\poset$ iff
$p:D\to 2$ for some finite $D\su\om\times\om$.

Each bijection $\tlpi:\om\to\om$ induces
an automorphism $\pi:\poset\to\poset$ defined by:
Given $p:D\to 2$ then
$\pi(p):E\to 2$
where $E=\{(\tlpi(i),j)\st (i,j)\in D\}$
and $\pi(p)(\tlpi(i),j)=p(i,j)$
for each $(i,j)\in D$.

Let $\gg$ be the group of automorphisms of $\poset$ 
generated by $\{\pi_{i,j}\st i<j<\om\}$ where $\tlpi_{i,j}:\om\to\om$
is the bijection which swaps $i$ and $j$.

The normal filter $\ff$ is generated by the subgroups $\{H_n\st n<\om\}$
where $H_n=\{\pi\in\gg\st \tlpi\res n =id\}$.
For $G$ $\poset$-generic over $M$, we let $\nn$ with
$M\su\nn\su M[G]$ be the symmetric model determined
by $(G,\gg,\ff)$, so $M\su\nn\su M[G]$. The model $\nn$ is the 
Basic Cohen model for the negation of the axiom of choice.
In $M[G]$ we define
$$x_n=\{k<\om\st \exists p\in G\;\; p(n,k)=1\}
\rmand
X=\{x_n\st n<\om\}.$$

The set $X$ is in $\nn$ and $\nn$ thinks it is
Dedekind finite, so no enumeration of it is there.
Recall that in $\nn$ we define the poset $\poset_0=\inj(X,X)$ 
to be the set of all finite partial one-to-one maps from $X$ to $X$.
If $G_0$ is $\poset_0$-generic over $\nn$, then $\bigcup G_0$
will be the graph of a bijection from $X$ to $X$.

In both posets $\poset$ and $\poset_0$
the trivial condition is the empty set, i.e.,
$\one=\emptyset$
and a universal name for the empty set is also the empty set.
The standard names for elements of the ground model
are defined by induction as
$\check{x}=\{(\one,\check{y})\st y\in x\}$.
The names for unordered and ordered pairs are
$$\{\tau_1,\tau_2\}^\circ=\{(\one,\tau_1),(\one,\tau_2)\}
\rmand
(\tau_1,\tau_2)^\circ=\{(\one,\{\tau_1\}^\circ),
(\one,\{\tau_1,\tau_2\}^\circ)\}.$$

\bigskip
Working in $\nn$ let
$$\Gamma=\{(r,\check{r})\st r\in\poset_0\}$$
be the usual name for $G_0$, the $\poset_0$-generic filter over $\nn$.

\bigskip
Working in $M$
let $\name{\Gamma}$ be a hereditarily symmetric
$\poset$-name\footnote{Yes, that's right, the name of a name.}
for $\Gamma$.
Let $\name{\poset_0}$ be a hereditarily symmetric name
for $\poset_0$.
Let
$$\name{x}_n=\{(p,\check{k})\st p\in\poset\rmand p(n,k)=1\}.$$

For each $n$ let
$$\namecheck{x}_n=\{(p, (\one,\checkcheck{k})^\circ)\st p(n,k)=1\}.$$
This will be a $\poset$-name for $\check{x}_n$ the standard
$\poset_0$-name for $x_n$.
This means that if $G$ is $\poset$-generic
over $M$ then $\namecheck{x}_n^G=\check{x}_n$, i.e, the standard
name of $x_n$ not $x_n$.
Note that if $\tlpi$ maps column $m$ to column $m^\pr$, then
$$ \pi(\name{x}_m)=\name{x}_{m^\pr} \;\;\rmand\;\;
\pi(\namecheck{x}_m)=\namecheck{x}_{m^\pr}.$$
For $\si\in\inj(\om,\om)$ 
(the graph of a finite injection)
define
$$\name{r}_\si=\{(\one,(\name{x}_i,\name{x}_j)^\circ)\st (i,j)\in\si\}.$$
Note that for any $p\in\poset$ and $\poset$-name $r$ if
$p\forces r\in \name{\poset}_0$, then there exists $q\leq p$ and 
$\si\in\inj(\om,\om)$
 such that $q\forces r=\name{r}_\si$.

\bigskip
Back working in $\nn$ note that
$$\one\forces_{\poset_0}\exists u \;\;(\exists v \;\;u\neq v\rmand
\exists r\in\Gamma \;\;(u,v)\in r)$$
write this as
$$\one\forces_{\poset_0}\exists u \;\;\theta(u,\Gamma).$$
We claim that there does not exists a $\poset_0$-name $\tau$ in $\nn$
such that
$$\nn\models\mbox{``}\;\one\forces_{\poset_0}\theta(\tau,\Gamma)\;\mbox{''}$$
and hence the maximal principle fails.
Suppose not and let $\name{\tau}$ be a
hereditarily symmetric $\poset$-name
for $\tau$.

\bigskip
Take $p\in G$ such that
$$p\forces \name{\nn}\models 1\forces_{\name{\poset_0}}
\theta(\name{\tau},\name{\Gamma}).$$
Working in $M$
choose $n$ so that $\dom(p)\su n\times\om$ and for every $\pi\in H_n$
$\;\;\pi(\name{\tau})=\name{\tau}$ and
$\pi(\name{\poset_0})=\name{\poset_0}$.

\bigskip
Working in $\nn$ let
$r_{id_n}=\{(x_i,x_i)\st i<n\}$.
We can find $r\leq r_{id_n}$ and $\check{x}_m$ such that
$$r\forces \tau=\check{x}_m.$$
Note that $\nn $ will not know which subscript goes with which
element of $X$ but we know
that $m\geq n$.

\bigskip
Working back in $M$ find $q\leq p$ and $\si\in \inj(\om,\om)$
with $\si\supseteq id_n$ such that
$$q\forces \name{\nn}\models \;\name{r}_\si\forces_{\name{\poset}_0}
\name{\tau}=\namecheck{x}_m$$
We write this as:
$$q\forces \psi(\name{\nn},\name{r}_\si,\name{\poset}_0,
\name{\tau},\namecheck{x}_m)$$
Now take $N>n$ with $dom(q)\su N\times\om$, $n\leq m<N$,
and $\si\su N\times N$.   Let $\pi\in\gg$ be determined
by the bijection $\tlpi:\om \to \om$
given by swapping the interval of columns $[n,N)$ with $[n+N,2N)$, 
i.e., swap
$k$ and $N+k$ for each $k$ with $n\leq k<N$.
Note that the corresponding automorphism $\pi$ of
$\poset$ has the property $\pi(\namecheck{x}_m)=\namecheck{x}_{m+N}$.
Let
$$\si^\pr=id_n\cup
\{(i+N,j+N)\st (i,j)\in\si \rmand i,j\geq n\}$$
and note that $\pi(\name{r}_\si)=\name{r}_{\si^\pr}$.
Since $\pi\in H_n$ it fixes $\name{\poset}_0$
and $\name{\tau}$ so
$$\pi(q)\forces \psi(\name{\nn},\name{r}_{\si^\pr},\name{\poset}_0,
\name{\tau},\namecheck{x}_{m+N}).$$

\bigskip
But $q$ and $\pi(q)$ are compatible so we may find $G$ which
is $\poset$-generic over $M$ containing them both.  In the model
corresponding model $\nn$ we will get that
$$r_\si\forces \tau=\check{x}_m\;\;\;\;\rmand\;\;\;\; 
r_{\si^\pr}\forces \tau=\check{x}_{m+N}$$
but this is a contradiction because $r_\si$ and $r_{\si^\pr}$ are
compatible.
\qed

\begin{example}
Recall that $Fn(I,J,\om)$ is the partial 
order of finite maps from $I$ to $J$, i.e. $r\su I\times J$ is
finite and $(u,v)\in r$ and $(u,w)\in r$ implies $v=w$.
Some other posets in $\nn$ for which the maximum principle fails
and for which some variant of the above argument works are:
\begin{enumerate}
\item $Fn(X,2,\om)$ $\;\;\;\;\;\;$ $\exists u \;\;(
\exists r\in\Gamma \;\;(u,0)\in r)$
\item $Fn(X,\om,\om)$ $\;\;\;\;\;\;$ $\exists u \;\;(
\exists r\in\Gamma \;\;(u,0)\in r)$
\item $Fn(X,X,\om)$
$\;\;\;\;$ $\exists u \;\;(\exists r\in\Gamma \;\;(u,x_0)\in r)$
\end{enumerate}
\end{example}
Proofs are left for the reader.
Finally we show that in $\nn$ the maximum principle holds for 
$\poset_1=\inj(\om,X)$.  Recall that this is the partial order
of the finite one-to-one maps from
$\om$ into $X$.
The key to the proof is Lemma \ref{key}, but
first we note some preliminary lemmas.

Define $H_n^\infty$ to be the subgroup of automorphisms of $\poset$ which
are determined by bijections $\tlpi:\om\to\om$ which are
the identity on $n$ ,i.e., $\tlpi(i)=i$ for all $i<n$.
Hence $H_n$ is $\gg\cap H_n^\infty$.
The elements of $H_n^\infty$ do not have to be in the ground model $M$
or even $M[G]$.

\begin{lemma}\label{decompose}
Suppose $k>n$ and $\pi\in H_n^\infty$
then there exists $\pi_1\in H_n$ and $\pi_2\in H_k^\infty$ such that
$\pi=\pi_1\circ\pi_2$.
\end{lemma}
\proof
Consider any orbit of $\tlpi$
which contains at least one of the $j<k$.
If it is finite, we set $\tlpi_1=\tlpi$ on it and
put $\tlpi_2$ to be the identity.
If it is an infinite orbit, write it as
$\{a_m\st m\in \zz\}$ where $\tlpi(a_m)=a_{m+1}$.
Since there are only finitely many $a_i$ with $0\leq a_i<k$,
we may renumber them so that for some $N$
any $a_i$ with $0\leq a_i<k$ is
in the set $a_1,\ldots, a_{N-1}$.
On this orbit define
$\tlpi_1$ to shift the list $a_1,a_2,\ldots, a_{N}$ up
one and send the last to the beginning, i.e.,
$\tlpi_1(a_i)=a_{i+1}$
for $1\leq i<N$ and $\tlpi_1(a_{N})=a_1$.
Define $\tlpi_2$ to shift the $\zz$-chain:
$$\ldots,a_{-2},a_{-1},a_{0},a_{N},a_{N+1},\ldots$$ i.e.,
$\tlpi_2(a_j)=a_{j+1}$ except when $j=0$ and
then $\tlpi_2(a_0)=a_{N}$.
\qed

\begin{lemma}\label{symmetricname}
For any hereditarily symmetric $\poset$-name
$\tau$, if every $\pi\in H_n$ fixes $\tau$, i.e., $\pi(\tau)=\tau$,
then every
$\pi\in H_n^\infty$ fixes $\tau$.
\end{lemma}
\proof
This is proved by induction on the rank of $\tau$.  Suppose
that $\pi\in H_n^\infty$ and $(p,\si)\in\tau$.  Choose $k>n$ so
that $\dom(p)\su k\times\om$ and $H_k$ fixes $\si$.  By Lemma
\ref{decompose}
there exists $\pi_1\in H_n$ and $\pi_2\in H_k^\infty$ such that
$\pi=\pi_1\circ\pi_2$.  It follows that
$(\pi(p),\pi(\si))=(\pi_1(p),\pi_1(\si))$ since $\tlpi_2$
is that identity on $k$, so $\pi_2(p)=p$, and since by induction on
rank $\pi_2(\si)=\si$.  Since $\pi_1$ fixes $\tau$ we have that
$(\pi(p),\pi(\si))\in\tau$.  It follows that
$\pi(\tau)\su\tau$.  Applying the same argument to $\pi^{-1}$ shows
that $\pi^{-1}(\tau)\su\tau$ and therefore $\tau\su\pi(\tau)$ and
so $\pi(\tau)=\tau$.
\qed

\begin{lemma}\label{key}
Suppose $G$ is $\poset$-generic over $M$ and $\nn=\nn_G$ is
the symmetric inner model with $M\su\nn\su M[G]$.
Working in $M[G]$ define $$x_i=\{j\in\om\st \exists p\in G\; p(i,j)=1\}$$
and let
$$G_1=\{r\in\poset_1\st \forall i\in \dom(r)\;\; r(i)=x_i\}.$$
Then $G_1$ is $\poset_1$-generic over $\nn$ and $\nn[G_1]=M[G]$.

\medskip
\noindent Conversely, if $\tlG_1$ is $\poset_1$-generic over
$\nn$, then
$$\tlG=\{s\in\poset\st \forall (i,j)\in\dom(s)
\;\;[s(i,j)=1 \rmiff \exists p\in \tlG_1
\; j\in p(i)]\}$$
is $\poset$-generic over $M$ and
$\nn=\nn_{\tlG}$.
\end{lemma}
\proof
First we see that $G_1$ is $\poset_1$-generic
over $\nn$.  In this proof we will use
$r_\si\in \poset_1$ for $\si\in\inj(\om,\om)$ to
refer to the condition satisfying
$r_\si(i)=x_{\si(i)}$ for each $i\in\dom(\si)$.

Working in $M$ suppose that $\name{D}$ is a symmetric
name and $s\in \poset$ satisfies:
$$s\forces \name{D}\;\su\;\name{\poset}_1\mbox{ is dense open.}$$
Choose $n$ so that every $\pi$ in $H_n$ fixes $\name{D}$ and 
$\dom(s)\su n\times\om$.
Choose $t\tleq s$, $m>n$, and a one-to-one $\si:m\to\om$ such that
$\si\supseteq id_n$ and
$$t\forces \name{r}_\si\in \name{D}$$
where 
$$\name{r}_\si=\{(\check{j},\name{x}_{\si(j)})^\circ\st j<m\}.$$
Let $\pi\in H_n$ be an automorphism for which
$\tlpi(\si(j))=j$ for every $j<m$.  It follows that
$$\pi(\name{r}_\si)=\name{r}_{id_m}$$
and
$$\pi(t)\forces \name{r}_{id_m}\in \name{D}.$$
Since $\pi(t)\tleq s$ and $s$ and $D$ were arbitrary it
follows that $G_1$ meets every dense subset of $\poset_1$ in $\nn$.

Since $M[G]$ is the smallest model of ZF containing  $G$ and
including $M$ we have that $M[G]\su\nn[G_1]$.  The other
inclusion follows since $G_1$ is easily definable from $G$.

Next we prove the ``Conversely'' statement.  Suppose that
$D\su\poset$ is dense and in $M$.  We must show it meets $\tlG$.

Working in $\nn$  for $s\in\poset$ and
$q\in\poset_1$ define $s\sqleq q$ as follows: 
For any $(i,j)\in \dom(s)$ we
have that $i\in\dom(q)$ and ($s(i,j)=1$ iff $j\in q(i)$).

We claim that
$$E=\{q\in\poset_1\st\exists s\in D\;\; s\sqleq q\}$$
is dense in $\poset_1$.  Since $E$ is in $\nn$ we have
that $E$ meets $\tlG_1$.  It follows that $D$ meets $\tlG$.

To prove $E$ is dense work in $M[G]$.  Fix $p\in\poset_1$.
Take $\pi\in\gg$
so that $p(i)=x_{\tlpi(i)}$ for each $i\in\dom(p)$.
Since $D$ is dense, so is $\pi^{-1}(D)$. Take
$s\in G\cap \pi^{-1}(D)$.  Then $\pi(s)\in D$ and
if $s_0=s\res\dom(p)$, then $s_0\sqleq p$.  By genericity it
is easy to find $q\tleq p$ with $\pi(s)\sqleq q$.

Finally, we show $\nn_G=\nn_{\tlG}$.  Let $\tlpi:\om\to\om$
be the bijection defined by $\tlpi(i)=j$ iff
$\exists p\in G_1$ with $p(i)=x_j$.   Then $\pi\in H^\infty_0$.
Note also that $\tlG=\pi(G)$.

It is a standard fact that the hereditarily
symmetric $\poset$-names in $M$ are closed under
$\gg$.
Combining Lemmas \ref{decompose} and \ref{symmetricname} gives
that the same is true for any $\pi\in H^\infty_0$.
To see this, suppose $\tau$ is fixed by $H_n$.  Decompose
$\pi=\pi_1\circ\pi_2$ with $\pi_2\in H_n^\infty$ and $\pi_1\in \gg$.
Then $\pi(\tau)=\pi_1(\tau)$.

Note that we have that
$$\tau^G=\pi(\tau)^{\pi(G)}=\pi_1(\tau)^{\tlG}$$
and hence $\nn_G\su\nn_{\tlG}$.  Similarly $\nn_{\tlG}\su\nn_{G}$
so they are equal.
\qed

\begin{theorem}
In $\nn$ the partial order $\poset_1=\inj(\om,X)$ satisfies 
the maximum principle.
\end{theorem}
\proof

Let $\pnames$ be the class\footnote{
This may be assumed to be a definable class in 
$M[G]$ and in $\nn$.
It is easy to see this would be true if we make the additional
assumption that $M$ is a model of $V=L$.  In
general one can make it true by adding 
a unary predicate for $M$
to the models.
See Solovay \cite{solovay} p.5-6.}
of $\poset$-names in $M$.

Working in $\nn$ define a mapping which
takes $\pnames$ to $\poset_1$-names as follows:
$$\hat{\tau}=\{(q,\hat{\si})\st
\exists r\; (r,\si)\in\tau \rmand r\sqleq q\}.$$
The relation $\sqleq$ is defined in the proof of Lemma \ref{key}.
It then follows that
$$\hat{\tau}^{\tlG_1}={\tau}^{\tlG}$$
for any $\tlG_1$ which is $\poset_1$-generic
over $\nn$ and $\tlG$ defined from it as in Lemma \ref{key}.

In $\nn$ suppose that
$$p_0\forces_{\poset_1}\;\exists x \;\theta(x).$$
For any $\tlG_1$ $\poset_1$-generic over $\nn$ with $p_0$ in $\tlG_1$,
we know that
$$\nn[\tlG_1]\models \exists x \;\theta(x)$$
by the definition of forcing.
By the key Lemma \ref{key}, $\nn[\tlG_1]=M[\tlG]$ and
so for some $\tau$ in $\pnames$
$$M[\tlG]\models \theta(\tau^{\tlG})$$
and so
$$\nn[\tlG_1]\models \;\theta(\hat{\tau}^{\tlG_1}).$$

It follows that in $\nn$
$$\forall q\tleq p_0\;\exists r\tleq q\;\exists \tau\in \pnames\;\;
r\forces_{\poset_1}\;\theta(\hat{\tau}).$$
By using the replacement axiom in $\nn$ and the axiom of choice in
$M$ we can find $\la\tau_\al\st\al<\ka\ra\in M\su\nn$ such that
in $\nn$:
$$\forall q\tleq p_0\;\exists r\tleq q\;\exists \al<\ka\;\;
r\forces_{\poset_1}\;\theta(\hat{\tau}_\al).$$
But this existential quantifier is essentially over an ordinal,
so by a proof similar to Proposition \ref{wellorder} we
can find a name $\tau$ such that
$$p_0\forces\theta(\tau)$$
and the maximum principle is proved.

Working in $\nn$ the name $\tau$ can be found as follows.
Let
$$\rho=\{(q,\hat{\tau}_\al)\st q\tleq p_0,\; 
q\forces\theta(\hat{\tau}_\al), \rmand \forall \be<\al\;\;
q\forces\;\neg\theta(\hat{\tau}_\be)\}.$$
Then $\rho$ is the name of a singleton $\{u\}$ where $u$ satisfies
$\theta$.  As in the usual proof of the maximum principle,
to remove the enclosing braces note that $u=\cup\{u\}$, so letting
$$\tau=\cup^\circ\rho=
\{(q_3,\si_2)\st \exists (q_1,\si_1)\in\rho \;\;\exists
q_2\; (q_2,\si_2)\in \si_1 \;\; q_3\tleq q_1,q_2\}$$
does the job.

\qed

\address

\end{document}